\documentclass{amsart}
\newtheorem{Thm}{Theorem}[section]
\newtheorem{Prop}[Thm]{Proposition}
\newtheorem{Cor}[Thm]{Corollary}
\newtheorem{Lem}[Thm]{Lemma}
\numberwithin{equation}{section}
\begin{document}

\title[$p$-harmonic boundary for graphs and manifolds]
{The $p$-harmonic boundary for quasi-isometric graphs and manifolds}

\author[M. J. Puls]{Michael J. Puls}
\address{Department of Mathematics \\
John-Jay College-CUNY \\ 
445 West 59th Street \\
New York, NY 10019, USA}
\email{mpuls@jjay.cuny.edu}

\begin{abstract} 
Let $p$ be a real number greater than one. Suppose that a graph $G$ of bounded degree is quasi-isometric with a Riemannian manifold $M$ with certain properties. Under these conditions we will show that the $p$-harmonic boundary of $G$ is homeomorphic to the $p$-harmonic boundary of $M$. We will also prove that there is a bijection between the $p$-harmonic functions on $G$ and the $p$-harmonic functions on $M$.
\end{abstract}

\date{February 3, 2010}
\maketitle
\section{Introduction}\label{Introduction}
Let $(X, d_X)$ and $(Y, d_Y)$ be metric spaces. A map $\phi \colon X \rightarrow Y$ is called a {\em quasi-isometry} if it satisfies the following two conditions:
\begin{enumerate} 
  \item There exists constants $a \geq 1, b \geq 0$ such that for $x_1, x_2 \in X$
                   \[ \frac{1}{a} d_X(x_1, x_2) - b \leq d_Y( \phi(x_1), \phi(x_2)) \leq a d_X(x_1, x_2) + b. \]
  \item There exists a positive constant $c$ such that for each $y \in Y$, there exists an $x \in X$ that satisfies $d_Y(\phi(x), y) < c.$
\end{enumerate} 
Let $G$ be a graph and let $x$ be a vertex of $G$. The set of neighbors of $x$ will be denoted by $N_x$ and $deg(x)$ will denote the number of neighbors of $x$. We shall say that $G$ is of {\em bounded degree} if there exists a positive integer $k$ such that $deg(x) \leq k$ for every vertex $x$ of $G$. A path in $G$ is a sequence of vertices $x_1, x_2, \dots, x_n$ where $x_{i+1} \in N_{x_i}$ for $1 \leq i \leq n-1$. A graph $G$ is connected if any two given vertices of $G$ are joined by a path. All graphs considered in this paper will be countably infinite, connected, of bounded degree with no self-loops. Two vertices $x$ and $y$ in $G$ are connected by an edge if and only if $y \in N_x$. Assign length one to each edge of $G$, then $G$ is a metric space with respect to the shortest path metric. Let $d_G(\cdot, \cdot)$ denote this metric. So if $x$ and $y$ are vertices in $G$, then $d_G(x,y)$ is the length of the shortest path joining $x$ and $y$. We will drop the subscript $G$ from $d_G( \cdot, \cdot)$ when it is clear what graph $G$ we are working with. If $A$ is a set of vertices of $G$ then $\#A$ will denote the cardinality of $A$.

Let $M$ be a complete, connected, and non-compact, smooth Riemannian manifold of dimension $n \geq 2$. Using the Riemannian distance, which we will denote by $d_M(\cdot, \cdot)$, $M$ is also a metric space. We will use $dx$ for the Riemannian volume element. For $x \in M, B_r(x)$ will denote the metric ball centered at $x$ of radius $r$; and $Vol(S)$ will be the volume of a measurable set $S \subseteq M$. In addition to the conditions in the first sentence of this paragraph, all manifolds considered in this paper will also have the following properties:
\begin{description}
  \item[(V)] There are positive increasing functions $V_0(r)$ and $V_1(r)$ on $(0, \infty)$ that satisfy
    \[ V_0(r) \leq Vol(B_r(x)) \leq V_1(r) \]
for all $x \in M$.
  \item[(P)] For $r > 0$, there exists a real number $C_r$ such that for any $y \in M$ and any smooth function $f$ on $B_r(y)$
     \[ \int_{B_r(y)} \vert f(x) - \bar{f} \vert \, dx \leq C_r\int_{B_r(y)} \vert \nabla f(x) \vert \, dx, \]
where $\bar{f} = (Vol (B_r(y)))^{-1} \int_{B_r(y)} f(x) \, dx.$
\end{description} 
These properties are satisfied by any complete manifold $M$ where the Ricci curvature of $M$ is uniformly bounded from below by $-(n-1)K^2$, where $K > 0$, and the injective radius of $M$ is positive. 

Let $p$ be a real number greater than one. In section \ref{preliminaries} we will define the $p$-harmonic boundary for both graphs and manifolds. It was shown in \cite[Theorem 2.7]{PulsPJM} that if $G$ and $H$ are quasi-isometric graphs, then their $p$-harmonic boundaries are homeomorphic. On the other hand Theorem 2 of \cite{Lee} says that if $M$ and $N$ are quasi-isometric manifolds, then their $p$-harmonic boundaries are homeomorphic. A reasonable question to ask is the following: if a graph $G$ of bounded degree is quasi-isometric with a complete Riemannian manifold $M$, how are their $p$-harmonic boundaries related? Reinterpreting Theorem 2 of \cite{Kanai} into our setting it was shown that if $G$ is quasi-isometric with $M$, then the $p$-harmonic boundary of $G$ is empty if and only if the $p$-harmonic boundary of $M$ is empty. Theorem 1.1 of \cite{Holasor} says that if there is a quasi-isometry from $G$ to $M$ then the $p$-harmonic boundary of $G$ contains one element if and only if the $p$-harmonic boundary of $M$ contains one element. In this paper we extend these results by proving
\begin{Thm} \label{boundhomeo} 
Let $G$ be a graph of bounded degree and let $M$ be a complete Riemannian manifold with dimension at least two and that has properties (V) and (P). If $G$ and $M$ are quasi-isometric, then their $p$-harmonic boundaries are homeomorphic.
\end{Thm}
In Section \ref{preliminaries} we will also define what it means for a function to be $p$-harmonic on $G$ and $M$. Our other main result for this paper is:
\begin{Thm} \label{bijectionpharm}
Let $G$ be a graph of bounded degree and let $M$ be a complete Riemannian manifold with dimension at least two and that has properties (V) and (P). If $G$ and $M$ are quasi-isometric, then there is a bijection between the bounded $p$-harmonic functions on $G$ and the bounded $p$-harmonic functions on $M$.
\end{Thm}

This paper is organized as follows: In Section \ref{preliminaries} we give some preliminaries concerning the $p$-harmonic boundary and $p$-harmonic functions. In section \ref{nets} we define $\kappa$-nets and give some results concerning $\kappa$-nets that we will need. We prove our main results in Section \ref{proofsofmain}.

I would like to thank the referee for some useful remarks concerning this paper.

This work was partially supported by PSC-CUNY grants 60123-38 39 and 62598-00 40.

\section{Preliminaries} \label{preliminaries}
Let $1 < p \in \mathbb{R}$. In this section we will define the $p$-harmonic boundary and $p$-harmonic functions. Furthermore, we will state some properties of these concepts that will be needed later in the paper. We will also set some notation to be used in this paper, and give some facts that will be needed concerning estimates on volumes of metric balls. We begin by defining certain function spaces that will be used in our definitions. For more detailed explanations about these function spaces and about the $p$-harmonic boundary see Section 1 of \cite{PulsPJM} for graphs and for Riemannian manifolds see Section 1 of \cite{Lee} or Chapter III.1 of \cite{SarNak}.

Let $M$ be a complete Riemannian manifold. Denote by $D_p(M)$ the set of continuous real-valued functions on $M$ for which $\nabla f \in L^p(M)$, where $\nabla f$ is the distributional gradient of $f$. Set $BD_p(M)$ equal to the set of bounded functions in $D_p(M)$. Under the usual operations of function addition, pointwise multiplication of functions and scalar multiplication $BD_p(M)$ is a commutative algebra. Furthermore, $BD_p(M)$ is a Banach algebra with respect to the following norm
\[ \parallel f \parallel_{BD_p} = \parallel \nabla f \parallel_p + \parallel f \parallel_{\infty} \]
where $\parallel \cdot \parallel_{\infty}$ denotes the sup-norm and $\parallel \cdot \parallel_p$ is the $L^p$-norm. Let $C_c(M)$ be the set of continuous functions on $M$ with compact support. Denote by $B(\overline{C_c(M)}_{D_p})$ the closure of $C_c(M)$ in $BD_p(M)$ with respect to the following topology: A sequence $(f_n)$ in $C_c(M)$ converges to $f \in BD_p(M)$ if $\sup_K \mid f_n - f \mid \rightarrow 0$ as $n \rightarrow \infty$ for each compact subset $K$ in $M$, $(f_n)$ is uniformly bounded on $M$ and 
  \[ \lim_{n \rightarrow \infty} \int_M \mid \nabla (f_n - f)(x) \mid^p \, dx \rightarrow 0. \] 

We now proceed to define analogous function spaces for a graph $G$ of bounded degree. Let $V$ be the vertex set of $G$ and let $x \in V$. For a real-valued function $f$ on $V$ we define the $p$-th power of the gradient and the $p$-Dirichlet sum by 
\begin{equation*} 
  \begin{split}
  \vert Df(x) \vert^p &  =    \sum_{y \in N_x} \vert f(y) - f(x) \vert^p, \\
    I_p(f, V)         & =    \sum_{x \in V} \vert Df(x) \vert^p. 
  \end{split}
\end{equation*}
In this setting $D_p(G)$ will be the set of functions $f$ for which $I_p(f, V) < \infty$. Under the following norm $D_p(G)$ is a reflexive Banach space
\[ \parallel f \parallel_{D_p} = \left( I_p(f,V) + \vert f(o) \vert^p \right)^{1/p}, \]
where $o$ is a fixed vertex of $G$. Let $BD_p(G)$ be the set of bounded functions in $D_p(G)$. It is also the case that $BD_p(G)$ is a commutative algebra under the operations of function addition, pointwise multiplication and scalar multiplication. With respect to the following norm $BD_p(G)$ is a Banach algebra
\[ \parallel f \parallel_{BD_p} = \left( I_p(f,V) \right)^{1/p} + \parallel f \parallel_{\infty}, \]
where $\parallel \cdot \parallel_{\infty}$ is the usual sup-norm. The set $C_c(G)$ will consist of all functions on $V$ with compact support. Denote by $B(\overline{C_c(G)}_{D_p})$ the closure of $C_c(G)$ in $BD_p(G)$ with respect to the $D_p$-norm.

In what follows $X$ will be either $M$ or $G$. A character on $BD_p(X)$ is a nonzero homomorphism from $BD_p(X)$ into the complex numbers. Denote by $Sp( BD_p(X))$ the set of characters on $BD_p(X)$. With respect to the weak $\ast$-topology, $Sp(BD_p(X))$ is a compact Hausdorff space. The space $Sp (BD_p(X))$ is known as the spectrum of $BD_p(X)$. Let $C(Sp(BD_p(X)))$ denote the set of continuous functions on $Sp(BD_p(X))$. For each $f \in BD_p(X)$ a continuous functions $\hat{f}$ can be defined on $Sp(BD_p(X))$ by $\hat{f}(\tau) = \tau(f)$. Each $x \in X $ defines an element in $Sp(BD_p(X))$ via evaluation by $x$; that is, if $f \in BD_p(X)$, then $x(f) = f(x)$. It turns out that under this identification $X$ is an open dense subset of $Sp(BD_p(X))$. The compact Hausdorff space $Sp(BD_p(X))\setminus X$ is known as the {\em $p$-Royden boundary} of $X$, which we will denote by $R_p(X)$. Now, $B(\overline{C_c(X)}_{D_p})$ is closed in $BD_p(X)$ with respect to the $BD_p$-norm. The {\em $p$-harmonic boundary} of $X$ is the following subset of the $p$-Royden boundary
\[ \partial_p (X) \colon = \{ \tau \in R_p(X) \mid \hat{f}(\tau) = 0 \mbox{ for all } f \in B(\overline{C_c(X)}_{D_p})\}. \]
It can be shown that $\partial_p(X) = \emptyset$ if and only if $1 \in \overline{C_c(X)}_{D_p}$, where 1 is the constant function 1 on $X$. In this case it follows from  Theorem 2 of \cite{Kanai} that $\partial_p(G) = \partial_p(M)$ because $X$ is $p$-parabolic if and only if $1 \in \overline{C_c(X)}_{D_p}$. For the rest of this paper it will be assumed that $X$ is not $p$-parabolic. From now on we will implicitly assume the following
\begin{Lem} \label{gotobound}
Let $x \in \partial_p(X)$ and let $(x_n)$ be a sequence in $X$ that converges to $x$. Then $d_X(o, x_n) \rightarrow \infty$ as $n \rightarrow \infty$, where $o$ is a fixed point in $X$.
\end{Lem}
\begin{proof}
Suppose there exists a real number $M$ such that $d_X(o, x_n) \leq M$ for all $n \in \mathbb{N}$. Define a function $\chi_M \in C_c(X)$ by $\chi_M (y) = 1$ if $d_X(o,y) \leq M$ and $\chi_M(y) = 0$ if $d_X(o,y) > M$. Then $\hat{\chi}_M(x) = \lim_{n \rightarrow \infty} \chi_M (x_n) = 1$, a contradiction. Thus $d_X(o, x_n) \rightarrow \infty$ as $n \rightarrow \infty.$
\end{proof} 

Now suppose $X=M$ and let $W^{1,p}(M)$ be the set of functions $f \in L^p(M)$ for which $\nabla f \in L^p(M)$. If $h$ is a continuous function in $W^{1,p}_{loc}(M)$ that is a weak solution of 
\[ -\mbox{div}(\vert \nabla h \vert^{p-2} \nabla h ) = 0, \]
then we shall say that $h$ is $p$-harmonic.
On the other hand if it is the case $X=G$, then $h$ is defined to be $p$-harmonic if 
\[ \sum_{y \in N_x} \vert h(y) - h(x) \vert^{p-2}(h(y) - h(x)) = 0 \mbox{ for all }x \in V. \]
In the case $1 < p < 2$ we make the convention that $\vert h(y) - h(x) \vert^{p-2} (h(y) - h(x)) = 0$ if $h(y) - h(x) =0$. 

Let $BHD_p(X)$ be the set that consists of all bounded $p$-harmonic functions on $X$ that are contained in $D_p(X)$. We now state some properties of $p$-harmonic functions and the $p$-harmonic boundary that will be needed in the sequel. For proofs of these and other properties see Section 4 of \cite{PulsPJM} for graphs and Section 2 of \cite{Lee} for manifolds.
\begin{Thm} \label{decomp} ($p$-Royden decomposition) Let $f \in BD_p(X)$. Then there exists a unique $u \in B(\overline{C_c(X)}_{D_p})$ and a unique $h \in BHD_p(X)$ such that $f = u + h.$
\end{Thm}
Using the $p$-Royden decomposition the following characterization of functions in $BD_p(X)$ that vanish on $\partial_p(X)$ can be obtained.
\begin{Thm} \label{charvan}
Let $f \in BD_p(X)$. Then $f \in B(\overline{C_c(X)}_{D_p})$ if and only if $\hat{f}(\tau) = 0$ for all $\tau \in \partial_p(X).$
\end{Thm}
Observe that it follows immediately from the theorem that if $h \in BHD_p(X)$ is the $p$-harmonic function in the $p$-Royden decomposition of $f \in BD_p(X)$, then $f(\tau) = h(\tau)$ for all $\tau \in \partial_p(X)$. Furthermore, the following is also a consequence of the previous theorem 
\begin{Cor} \label{boundval} A function in $BHD_p(X)$ is uniquely determined by its values on $\partial_p(X)$. 
\end{Cor}
Note that if $\partial_p(X)$ contains only one element, then $BHD_p(X)$ consists precisely of the the constant functions on $X$.

\section{$\kappa$-nets}\label{nets}
For the proofs of our main results we will need to use $\kappa$-nets. In this section we will explain what a $\kappa$-net is, and give some of its properties that will be useful for our needs.

A {\em net} is a countable set $\Gamma$ with a family $\{N_g\}_{g \in \Gamma}$ of finite subsets $N_g$ of $\Gamma$ such that for all $g, h \in \Gamma, g \in N_h$ if and only if $h \in N_g$. For $g \in \Gamma$, each element of $N_g$ is called a neighbor of $g$. It is possible to think of a net as a countable graph with vertex set $\Gamma$ by connecting vertices $g$ and $h$ in $\Gamma$ by an edge if $h \in N_g$. Thus the definitions given in Section \ref{Introduction} for properties of a graph, such as a path, carry over to a net. 

Let $M$ be a complete Riemannian manifold. We shall say that a subset $\Gamma$ of $M$ is $\kappa$-separated for $\kappa > 0$ if $d_M(g, h) \geq \kappa$ whenever $g$ and $h$ are distinct points of $\Gamma$. Now assume that $\Gamma$ is a maximal $\kappa$-separated subset of $M$. By setting $N_g = \{ h \in \Gamma \mid 0 < d_M(g,h) \leq 3\kappa \}$ for each $g \in \Gamma$ we have a net structure on $\Gamma$. We define $\Gamma$ to be a {\em $\kappa$-net} if $\Gamma$ is a maximal $\kappa$-separated subset of $M$ with the net structure given above. It is easy to see that a $\kappa$-net in $M$ is connected due to our standing assumption that $M$ is connected. If $\Gamma$ is a $\kappa$-net, then it is also a countable, connected graph. Thus, $\Gamma$ is a metric space with respect to the shortest path metric. A couple of facts about $\kappa$-nets that we will need later are:
\begin{enumerate} 
  \item Let $\Gamma$ be a $\kappa$-net in $M$. Then for a given $r>0$, there exists a constant $C_r$ such that $\# (\Gamma \cap B_r(x)) \leq C_r$ for all $x \in M$. 
  \item Let $\Gamma$ be a $\kappa$-net in $M$. Then the inclusion map $\iota \colon \Gamma \rightarrow M$ is a quasi-isometry.
\end{enumerate}
The first fact can be proved by using the argument from Lemma 2.3 of \cite{Kanaicapp} since $Vol(B_R(x)) \leq \frac{V_1(R)}{V_0(r)} Vol (B_r(x)),$ where $0 < r < R < \infty$.  This fact shows that $\Gamma$ is a graph of bounded degree; and that for each $x \in M$ there exists at most $C_r$ elements $g$ in $\Gamma$ for which $B_r(g)$ contains $x$. The second fact is \cite[Lemma 2.6]{Kanaicapp}, where it was assumed that the Ricci curvature was bounded from below; it was also proven in Lemma 2.13 of \cite{HoloRMI} without any curvature assumptions on $M$.

Let $\kappa$ be a small positive number and let $\Gamma$ be a $\kappa$-net in a complete Riemannian manifold $M$. Let $1 < p \in \mathbb{R}$. The rest of this section is devoted to describing how to map functions from $BD_p(\Gamma)$ to $BD_p(M)$, and vice-versa. For each $g \in \Gamma$, pick a smooth function $\eta_g \in C_c(M)$ such that $0 \leq \eta_g \leq 1, \eta_g = 1$ on $B_{\kappa}(g), \eta_g = 0$ outside of $B_{\frac{3\kappa}{2}} (g)$, and that $\vert \nabla \eta_g \vert \leq c$, where $c$ is a constant that does not depend on $g$. For $x \in M$ define 
\[ \xi_g (x) = \frac{\eta_g(x)}{\sum_{ h \in \Gamma} \eta_h(x)} . \]
Now $\vert \nabla \xi_g \vert$ is uniformly bounded. Indeed, 
\begin{eqnarray*}
  \vert \nabla \xi_g \vert &  \leq  & \vert \nabla \eta_g \vert \left(\sum_{h \in \Gamma} \eta_h\right)^{-1} + \eta_g \sum_{h \in \Gamma} \vert \nabla \eta_h \vert \left(\sum_{h \in \Gamma} \eta_h \right)^{-2} \\
    & \leq  & \vert \nabla \eta_g \vert + \sum_{h \in \Gamma} \vert \nabla \eta_h \vert  \\
    & \leq   &(k + 2)c,
\end{eqnarray*}
where $k$ is a constant that satisfies $\#N_g \leq k$ for all $g \in \Gamma$. Let $\bar{f} \colon \Gamma \rightarrow \mathbb{R}$. Define a smooth function $f \colon M \rightarrow \mathbb{R}$ by
\begin{equation} 
f(x) = \sum_{g \in \Gamma} \bar{f}(g) \xi_g (x),   \label{graphtorm}
\end{equation}where $x \in M$. We are now ready to state and prove
\begin{Prop} \label{definefrm}
If $\bar{f} \in BD_p(\Gamma)$, then $ f\in BD_p(M).$ 
\end{Prop}
\begin{proof}
Let $g \in \Gamma$ and suppose $x \in B_{\kappa}(g)$. Now 
\begin{eqnarray*}
\nabla f(x)  & = &\sum_{h \in N_g \cup \{g\}} \bar{f}(h) \nabla \xi_h (x) \\
             & = &\sum_{h \in N_g} (\bar{f}(h) - \bar{f}(g)) \nabla \xi_h (x).
\end{eqnarray*}
The last equality is due to $\sum_{h \in N_g \cup \{g\}} \bar{f}(g) \xi_h(x) = \bar{f}(g)$ and $\sum_{h \in N_g \cup \{g\}} \nabla \xi_h (x) =0$. We now obtain
\begin{eqnarray*}
 \vert \nabla f(x) \vert^p & \leq & \left( \sum_{h \in N_g} \vert (\bar{f}(h) - \bar{f}(g)) \nabla \xi_h (x) \vert \right)^p \\
                         & \leq & \left( \left( \sum_{h \in N_g} \vert \bar{f}(h) - \bar{f}(g) \vert^p \right)^{1/p} \left( \sum_{h \in N_g} \vert \nabla \xi_h (x) \vert^q \right)^{1/q} \right)^p \\
  & \leq & \vert D\bar{f}(g) \vert^p (ck^{1/q})^p,
\end{eqnarray*}
where $\frac{1}{q} + \frac{1}{p} = 1$ and $k$ is a constant with $\#N_g \leq k$ for all $g \in \Gamma$. It now follows that,
\[ \int_M \vert \nabla f(x) \vert^p dx \leq \sum_{g \in \Gamma} \int_{B_{\kappa}(g)} \vert \nabla f(x) \vert^p dx \leq (ck^{1/q})^p V_1(\kappa) \sum_{g \in \Gamma} \vert D\bar{f}(g)\vert^p \]
Hence, $f \in BD_p(M)$.
\end{proof}
\begin{Cor} \label{compactfrm}
If $\bar{f} \in B(\overline{C_c(\Gamma)}_{D_p})$, then $f \in B(\overline{C_c(M)}_{D_p}).$
\end{Cor} 
\begin{proof}
By the proposition $f \in BD_p(M)$. Now let $(\bar{f}_n)$ be a sequence in $C_c(\Gamma)$ that converges to $\bar{f}$. For each $n, f_n \in C_c(M)$ because $\bar{f}_n \in C_c( \Gamma)$. We will now show that $f \in B(\overline{C_c(M)}_{D_p})$. By using the argument from the above proposition we see that 
\[ \int_M \vert \nabla (f_n -f)(x) \vert^p dx \leq C I_p(\bar{f}_n - \bar{f}, \Gamma), \]
where $C$ is a constant. Consequently, $\int_M \mid \nabla (f_n - f)(x) \mid^p dx \rightarrow 0$ as $n \rightarrow \infty$. Let $K$ be a compact subset in $M$. Set $\Gamma_K = \{ g \in \Gamma \mid d_M(g, K) < \frac{3\kappa}{2} \}$. Now $\# \Gamma_K$ is finite because $K$ is compact. Let $x \in K$ and let $\epsilon > 0$. Since $\xi_g (x) = 0$ if $g \notin \Gamma_K$ we obtain
\[ \vert f_n (x) - f(x) \vert \leq \sum_{g \in \Gamma_K} \vert \bar{f}_n(g) - \bar{f}(g) \vert. \]
For $h \in D_p(\Gamma)$ and $g \in \Gamma$, there exists a constant $C_g$ depending on $g$ such that $\vert h(g) \vert \leq C_g \parallel h \parallel_{D_p}$. Thus $(\bar{f}_n) \rightarrow \bar{f}$ pointwise. Hence for each $g \in \Gamma_K$ there exists a number $N(g)$ such that for $n > N(g), \vert \bar{f}_n (g) - \bar{f}(g) \vert < \frac{\epsilon}{\# \Gamma_K}$. Set $N = \max_{g \in \Gamma_K} \{ N(g)\}$. So for $n > N, \vert f_n(x) - f(x) \vert < \epsilon$ for all $x \in K$. Thus $\sup_K \vert f_n - f \vert \rightarrow 0$ as $n \rightarrow \infty$ for each compact subset $K$ in $M$. By making slight modifications to the proof of Theorem 1G from page 153 of \cite{SarNak} it follows that $f \in B(\overline{C_c(M)}_{D_p})$.
\end{proof}

Let $f \in BD_p(M)$. Define a function $f^{\ast} \colon \Gamma \rightarrow \mathbb{R}$ by 
\begin{equation} 
f^{\ast} (g) = \frac{1}{Vol (B_{4\kappa}(g))} \int_{B_{4\kappa}(g)} f dx.   \label{rmtograph}
\end{equation}
We will now show $f^{\ast} \in BD_p(\Gamma)$. 
\begin{Prop} \label{definefgraph}
Let $f \in BD_p(M)$, then $f^{\ast} \in BD_p(\Gamma)$.
\end{Prop}
\begin{proof}
Let $g \in \Gamma$. By H\"{o}lder's inequality and property (P) we get the following
\begin{eqnarray*}
V_1(4\kappa)^{p-1} \int_{B_{4\kappa}(g)} \vert \nabla f(x) \vert^p dx & \geq & Vol(B_{4\kappa}(g))^{p-1} \int_{B_{4\kappa}(g)} \vert \nabla f(x)\vert^p dx \\
                  & \geq & \left( \int_{B_{4\kappa}(g)} \vert \nabla f(x) \vert dx \right)^p  \\
                  & \geq & C^p \left( \int_{B_{4\kappa}(g)} \vert f(x) - f^{\ast}(g) \vert dx \right)^p,
\end{eqnarray*}
where $C$ is a constant. Let $h \in N_g$, then $d_M(g,h) \leq 3\kappa$. Consequently both $B_{4\kappa}(g)$ and $B_{4\kappa}(h)$ are contained in $B_{7\kappa}(g)$. Letting $\beta = C^{-p}V_1(4\kappa)^{p-1}$ for convience we now obtain
\begin{eqnarray*}
2\beta\int_{B_{7\kappa}(g)} \vert \nabla f(x) \vert^p dx & \geq & \beta\left( \int_{B_{4\kappa}(g)}\vert \nabla f(x) \vert^p dx  + \int_{B_{4\kappa}(h)} \vert \nabla f(x) \vert^p dx \right) \\
  & \geq & \left( \int_{B_{4\kappa}(g)} \vert f(x) - f^{\ast}(g) \vert dx \right)^p \\ &  & \hspace{.4in} + \left( \int_{B_{4\kappa}(h)} \vert f(x) - f^{\ast}(h) \vert dx \right)^p \\
  & \geq & \frac{1}{2^{p-1}} \left( \int_{B_{4\kappa}(g)} \vert f(x) - f^{\ast}(g) \vert dx \right. \\ & & \hspace{.4in} \left. + \int_{B_{4\kappa}(h)} \vert f(x) - f^{\ast}(h) \vert dx \right)^p \\
  & \geq & \frac{1}{2^{p-1}} \left( \int_{B_{4\kappa}(g) \cap B_{4\kappa}(h)} \vert f^{\ast} (g) - f^{\ast} (h) \vert dx \right)^p \\
  & \geq & \frac{1}{2^{p-1}} V_0(\kappa)^p \left( \vert f^{\ast} (g) - f^{\ast} (h) \vert^p \right).
\end{eqnarray*}
The last inequality follows from $B_{\kappa}(g) \subseteq B_{4\kappa}(g) \cap B_{4\kappa}(h)$, and three inequalities up is Jensen's inequality. Due to $\Gamma$ having bounded degree there exists a constant $C_1$ that does not depend on $f$ or $g$ for which
\[ C_1 \int_{B_{7\kappa}(g)} \vert \nabla f(x) \vert^p dx \geq \sum_{h \in N_g} \vert f^{\ast}(g) - f^{\ast}(h) \vert^p. \]
Furthermore, we saw earlier that if $x \in M$ then there exists at most $C_{7\kappa}$ balls $B_{7\kappa}(g)$ that contain $x$, where $C_{7\kappa}$ is a constant that does not depend on $f$ or $g$. Hence,
\[ C_{7\kappa} \int_M \vert \nabla f(x) \vert^p dx \geq \sum_{g \in \Gamma} \int_{B_{7\kappa}(g)} \vert \nabla f(x) \vert^p dx. \]
Summing up we obtain
\[ C_2 \int_M \vert \nabla f(x) \vert^p dx \geq I_p (f^{\ast}, \Gamma), \]
where $C_2$ is a suitable constant. Therefore, $f^{\ast} \in BD_p(\Gamma)$. 
\end{proof}
\begin{Cor} \label{compactfgraph} If $f \in B( \overline{C_c(M)}_{D_p})$, then $f^{\ast} \in B(\overline{C_c(\Gamma)}_{D_p})$.
\end{Cor} 
\begin{proof} 
Let $f \in B(\overline{C_c(M)}_{D_p})$, then $f^{\ast} \in BD_p(\Gamma)$ by the proposition. We will now show that $f^{\ast}$ is also an element of $B(\overline{C_c(\Gamma)}_{D_p})$. Let $(f_n)$ be a sequence in $C_c(M)$ that converges to $f$. For each $n, f_n^{\ast} \in C_c(\Gamma)$ since $f_n$ has compact support and $\Gamma$ is $\kappa$-separated. Arguing as in the proposition it can be shown that there exists a constant $C$ such that 
\[ I_p(f_n^{\ast} - f^{\ast}, \Gamma) \leq C\int_M \vert \nabla(f_n - f)(x) \vert^p dx. \]
Thus $I_p( f_n^{\ast} - f^{\ast}, \Gamma) \rightarrow 0$ as $n \rightarrow \infty$. Let $o$ be a fixed vertex of $\Gamma$ and let $\epsilon > 0$. A calculation shows that 
\[ \vert f_n^{\ast}(o) - f^{\ast}(o) \vert \leq \frac{1}{V_0(4\kappa)} \int_{B_{4\kappa}(o)} \vert (f_n - f)(x) \vert dx. \]
Now $\vert f_n(x) - f(x) \vert < \epsilon$ for large $n$ and all $x \in B_{4\kappa}(o)$ due to the closure of $B_{4\kappa}(o)$ being compact in $M$. Hence $\vert f_n^{\ast} (o) - f^{\ast}(o) \vert^p \rightarrow 0$ as $n \rightarrow \infty$. Therefore, $f^{\ast} \in B(\overline{C_c(\Gamma)}_{D_p})$.
\end{proof}

\section{Proofs of Theorems \ref{boundhomeo} and \ref{bijectionpharm}}\label{proofsofresults}\label{proofsofmain}
Let $G$ be a graph of bounded degree and let $M$ be a Riemannian manifold. Let $\Gamma$ be a maximal $\kappa$-separated net in $M$, where $\kappa$ is a small positive number. Recall that $\Gamma$ can also be considered as graph of bounded degree with vertex set $\Gamma$. Assume for now that $\partial_p(\Gamma)$ is homeomorphic to $\partial_p(M)$. We saw in Section \ref{nets} that the embedding $\iota \colon \Gamma \rightarrow M$ is a quasi-isometry, so the graph $G$ is quasi-isometric with $\Gamma$ because the composition of quasi-isometries is a quasi-isometry. It follows from Theorem 2.7 of \cite{PulsPJM} that $\partial_p(G)$ is homeomorphic to $\partial_p(M)$, as desired. In order to complete the proof of Theorem \ref{boundhomeo} we need to show that $\partial_p(M)$ is homeomorphic to $\partial_p(\Gamma)$, which we now proceed to do. 

We begin with two crucial lemmas.
\begin{Lem}\label{infinityconv}
Let $(y_n)$ be a sequence in $M$ with $d_M(o, y_n) \rightarrow \infty$ as $n \rightarrow \infty$, where $o$ is a fixed point in $M$. For each $n \in \mathbb{N}$ let $x_n \in \Gamma$ that satisfies $d_M(x_n, y_n) < \kappa$. Let $1 < p \in \mathbb{R}$. If $f \in BD_p(M)$, then $\vert f^{\ast}(x_n) - f(y_n) \vert \rightarrow 0$ as $n \rightarrow \infty$, where $f^{\ast}$ is defined by (\ref{rmtograph}).
\end{Lem}
\begin{proof}
For each $n \in \mathbb{N}$ we have that 
\begin{eqnarray*}
\vert f^{\ast}(x_n) - f(y_n) \vert & = & \left\vert \frac{1}{Vol(B_{4\kappa}(x_n))}\int_{B_{4\kappa}(x_n)} f(x) dx - f(y_n) \right\vert \\
                                  & = & \left\vert \frac{1}{Vol(B_{4\kappa}(x_n))} \int_{B_{4\kappa}(x_n)} \left(f(x) - f(y_n)\right) dx \right\vert \\
                                  & \leq & \frac{1}{Vol(B_{4\kappa} (x_n))} \int_{B_{4\kappa}(x_n)} \vert f(x) - f(y_n) \vert dx \\
                                   & \leq & \left( \frac{1}{Vol(B_{4\kappa}(x_n))} \int_{B_{4\kappa}(x_n)} \vert f(x) - f(y_n) \vert^p dx \right)^{1/p}.
\end{eqnarray*}
Let $r_x = d_M(y_n,x)$ and let $\gamma_x \colon [0, \infty) \rightarrow M$ be a geodesic parameterized by arclength that satisfies $\gamma_x(0) = y_n$ and $\gamma_x(r_x) = x$. By independence of path  
\[ \int_0^{r_x} \nabla f(\gamma_x(t)) \cdot \gamma'_x(t) dt = f(x) - f(y_n). \]
For notational convenience set $\beta = \frac{1}{Vol(B_{4\kappa}(x_n))}$. Consequently, 
\begin{eqnarray*}
\left( \beta \int_{B_{4\kappa}(x_n)} \vert f(x) - f(y_n) \vert^p dx \right)^{\frac{1}{p}} & \leq & \left( \beta \int_{B_{4\kappa}(x_n)} \left( \int_0^{5\kappa} \vert \nabla f(\gamma_x(t)) \vert dt \right)^p dx \right)^{1/p} \\
  & \leq & \int_0^{5\kappa} \left( \beta \int_{B_{4\kappa}(x_n)} \vert \nabla f(\gamma_x(t)) \vert^p dx \right)^{1/p} dt \\
& \leq & 5\kappa\beta^{1/p} \left( \int_{B_{5\kappa}(y_n)} \vert \nabla f(x) \vert^p dx \right)^{1/p}.
\end{eqnarray*}
The second to last inequality is Minkowski's Integral Inequality. By condition (V) we have a constant $V_0(5\kappa)$ such that $V_0(5\kappa) \leq Vol(B_{5\kappa}(x))$ for all $x \in M$. Hence,
\[ \vert f^{\ast} (x_n) - f(y_n) \vert \leq 5 \kappa V_0(5\kappa)^{-1/p} \left( \int_{B_{5\kappa}(y_n)} \vert \nabla f(x) \vert^p dx \right)^{1/p}. \]
Now $\int_{B_{5\kappa}(y_n)} \vert \nabla f(x) \vert^p dx \rightarrow 0$ as $d_M(0, y_n) \rightarrow \infty$ because $f\in BD_p(M)$. Therefore, $\vert f^{\ast} (x_n) - f(y_n) \vert \rightarrow 0$ as $n \rightarrow \infty.$
\end{proof} 

\begin{Lem} \label{convinftyrm} Let $\bar{f} \in BD_p(\Gamma)$ and let $f \in BD_p(M)$ be defined by (\ref{graphtorm}). Let $x \in \partial_p(\Gamma)$ and let $(x_n)$ be a sequence in $\Gamma$ that converges to $x$. Then $\vert f(x_n) - \bar{f}(x_n) \vert \rightarrow 0$ as $n \rightarrow \infty$. 
\end{Lem}
\begin{proof}
Let $n \in \mathbb{N}$ and let $\Gamma_{x_n} = \{g \in \Gamma \mid g \in B_{\frac{3\kappa}{2}}(x_n) \}$. Now,
\begin{eqnarray*}
\vert f(x_n) - \bar{f}(x_n) \vert & = & \left\vert \sum_{g \in \Gamma} (\bar{f}(g) - \bar{f}(x_n)) \xi_g (x_n) \right\vert \\
                                 & \leq & \left( \sum_{g \in \Gamma} \vert \bar{f}(g) - \bar{f}(x_n) \vert^p \xi_g(x_n) \right)^{1/p} \\
                                & \leq & \left( \sum_{g \in \Gamma_{x_n}} \vert \bar{f}(g) - \bar{f}(x_n) \vert^p \right)^{1/p}.
\end{eqnarray*}
The last inequality follows from $\xi_g(x_n) = 0$ if $g \notin \Gamma_{x_n}$. If $g \in \Gamma_{x_n} \setminus \{x_n\}$ then $g \in N_{x_n}$, so it follows that $\sum_{g \in \Gamma_{x_n}} \vert \bar{f}(g) - \bar{f}(x_n) \vert^p \rightarrow 0$ as $d_{\Gamma}(o, x_n) \rightarrow 0$ due to $\bar{f} \in BD_p(\Gamma)$. Hence, $\vert f(x_n) -\bar{f}(x_n) \vert \rightarrow 0$ as $n \rightarrow \infty$.
\end{proof}

We will now proceed to define a function $\Phi \colon \partial_p(\Gamma) \rightarrow \partial_p(M).$ Let $x \in \partial_p(\Gamma)$ and let $(x_n)$ be a sequence in $\Gamma$ that converges to $x$. Since $Sp(BD_p(M))$ is a compact Hausdorff space we may assume, by passing to a subsequence if necessary, that the sequence $(x_n)$ converges to a unique element $y \in Sp(BD_p(M))$. Define $\Phi(x) = y$. We now prove that $\Phi$ is well-defined and $y \in \partial_p(M).$
\begin{Prop} \label{welldefined} 
The map $\Phi$ is well-defined from $\partial_p(\Gamma)$ to $\partial_p(M)$.
\end{Prop}
\begin{proof}
Let $x$ and $y$ be as above. We first show that $\Phi$ is well-defined. Let $(x_n)$ and $(x'_n)$ be sequences in $\Gamma$ such that both $(x_n)$ and $(x'_n)$ converge to $x$. Suppose that $(x_n) \rightarrow y_1$ and $(x'_n) \rightarrow y_2$ in $Sp(BD_p(M))$ and further assume that $y_1 \neq y_2$. Let $f \in BD_p(M)$ that satisfies $f(y_1) =0$ and $f(y_2) = 1$. Define $f^{\ast}$ as in (\ref{rmtograph}). Setting $y_n = x_n$ in Lemma \ref{infinityconv} we obtain $\lim_{n \rightarrow \infty} f^{\ast}(x_n) = 0$ and $\lim_{n \rightarrow \infty} f^{\ast} (x'_n) = 1$. But $f^{\ast} \in BD_p(\Gamma)$ which means
\[ \lim_{n \rightarrow \infty} f^{\ast} (x_n) = f^{\ast}(x)= \lim_{n \rightarrow \infty} f^{\ast}(x'_n), \]
a contradiction. Thus $\Phi$ is well-defined.

We will now show that $\Phi(x) = y \in \partial_p(M)$. Suppose $y \notin \partial_p(M)$, then there exists an $f \in B(\overline{C_c(M)}_{D_p})$ such that $f(y) \neq 0$. Assume $f(y) > 0$ and define $f^{\ast}$ as in (\ref{rmtograph}). Since $f^{\ast} \in B(\overline{C_c(\Gamma)}_{D_p})$ it must be the case $\lim_{n \rightarrow \infty} f^{\ast}(x_n) = f^{\ast}(x) = 0$. However, Lemma \ref{infinityconv} says that $\lim_{n \rightarrow \infty} f^{\ast}(x_n) = f(y) > 0$, a contradiction. Thus $y \in \partial_p(M).$
\end{proof}

We now show that $\Phi$ is a bijection.
\begin{Prop} \label{bijection}
The map $\Phi$ is a bijection.
\end{Prop}
\begin{proof} 
We will begin by showing that $\Phi$ is one-to-one. Let $x_1, x_2 \in \partial_p(\Gamma)$ that satisfy $\Phi(x_1) = \Phi(x_2)$. Assume $x_1 \neq x_2$. Then there exits an $\bar{f} \in BD_p(\Gamma)$ with $\bar{f}(x_1) = 1$ and $\bar{f}(x_2) = 0$. Let $(x_n)$ and $(x'_n)$ be sequences in $\Gamma$ such that $(x_n) \rightarrow x_1$ and $(x'_n) \rightarrow x_2$. Using $\bar{f}$ define a function $f \in BD_p(M)$ by (\ref{graphtorm}). By Lemma \ref{convinftyrm} we see that $\lim_{n \rightarrow \infty} f(x_n) = 1$ and $\lim_{n \rightarrow \infty}f(x'_n) = 0$. This contradicts the assumption
\[ \lim_{n \rightarrow \infty} (x_n) = \Phi(x_1) = \Phi(x_2) = \lim_{n \rightarrow \infty} (x'_n). \]
Thus, $\Phi$ is one-to-one.

We will now show that $\Phi$ is onto. Let $y \in \partial_p(M)$ and let $(y_n)$ be a sequence in $M$ with $(y_n) \rightarrow y$. For each $n \in \mathbb{N}$, choose $x_n \in \Gamma$ that satisfies $d(x_n, y_n) < \kappa$. We claim that $(x_n) \rightarrow y$ in $Sp(BD_p(M))$. To see the claim let $f\in BD_p(M)$ and define $f^{\ast} \in BD_p(\Gamma)$ by (\ref{rmtograph}). By Lemma \ref{infinityconv} both $\vert f(x_n) - f^{\ast}(x_n) \vert \rightarrow 0$ and $\vert f^{\ast} (x_n) - f(y_n) \vert \rightarrow 0$ as $n \rightarrow \infty$. Hence 
\[ \vert f(x_n) - f(y_n) \vert \rightarrow 0 \mbox{ as } n \rightarrow \infty \]
for all $f \in BD_p(M)$. Thus $\lim_{n \rightarrow \infty} (x_n) = \lim_{n \rightarrow \infty} (y_n) = y$ and the claim is proved.

By passing to a subsequence if need be, we assume that the sequence $(x_n)$ converges to an unique element $x$ in the compact Hausdorff space $Sp(BD_p(\Gamma))$. To finish the proof we need to show $x \in \partial_p(\Gamma)$. Suppose $x \notin \partial_p(\Gamma)$, then $\bar{f}(x) \neq 0$ for some $\bar{f} \in B(\overline{C_c(\Gamma)}_{D_p})$. Using $\bar{f}$ define a function $f\in B(\overline{C_c(M)}_{D_p})$ via (\ref{graphtorm}). By Lemma \ref{convinftyrm}, $\lim_{n \rightarrow \infty} f(x_n) = \lim_{n \rightarrow \infty} \bar{f}(x_n)$, which implies $f(y) \neq 0$, contradicting $y \in \partial_p(M)$ and $f \in B(\overline{C_c(M)}_{D_p})$. Thus $x \in \partial_p(\Gamma)$ and $\Phi(x) = y$, which shows that $\Phi$ is onto.
\end{proof}

To finish the proof that the bijection $\Phi$ is a homeomorphism we only need to show that $\Phi$ is continuous, since both $\partial_p(\Gamma)$ and $\partial_p(M)$ are compact Hausdorff spaces. Let $U$ be an open subset of $\partial_p(M)$ and let $x \in \Phi^{-1}(U)$. Fix $\epsilon$ that satisfies $0 < \epsilon < 1$. Let $(x_n)$ be a sequence in $\Gamma$ for which $(x_n) \rightarrow x$ and let $y = \Phi(x) \in U$. By Proposition 1 of \cite{Lee} there exist a subset $\Omega$ of $M$ such that $y \in \overline{\Omega}$, where the closure is in $Sp(BD_p(M))$, and $\overline{\Omega} \cap \partial_p(M) \subseteq U$. It was shown in the proof of this proposition that $\Omega = \{ m \mid h(m) > \epsilon \}$, where $h \in BHD_p(M), 0 \leq h \leq 1, h(y) = 1$ and $h = 0$ on $\partial_p(M) \setminus U$. Using this $h$ a function $h^{\ast} \in BD_p(\Gamma)$ can be defined by (\ref{rmtograph}). Let $V = \{ z \in \partial_p(\Gamma) \mid h^{\ast}(z) > \frac{\epsilon}{2}\}$. The set $V$ is open due to $h^{\ast}$ being continuous. Moreover, it follows from Lemma \ref{infinityconv} that $x \in V$. Now let $z \in V$ and suppose $\Phi(z) \notin U$. Consequently, $h(z_n) \rightarrow 0$ where $(z_n)$ is a sequence in $\Gamma$ that converges to $z$. 
It follows from Lemma \ref{infinityconv} that $h^{\ast}(z_n) < \frac{\epsilon}{2}$ for large $n$, contradicting $z \in V$. Thus $V \subseteq \Phi^{-1}(U)$ and the continuity of $\Phi$ is established. Therefore, $\Phi \colon \partial_p (\Gamma) \rightarrow \partial_p(M)$ is a homeomorphism and Theorem \ref{boundhomeo} is proved.

The proof of Theorem \ref{bijectionpharm} pretty much follows the same path as the proof of Theorem \ref{boundhomeo}. Let $G$ be a graph of bounded degree and let $M$ be a complete Riemannian manifold. Let $\Gamma$ be a maximal $\kappa$-separated net in $M$, where $\kappa$ is a small positive constant. The graph with bounded degree $\Gamma$ is quasi-isometric with $G$ so by Theorem \cite[Theorem 2.8]{PulsPJM} there is a bijection between $BHD_p(G)$ and $BHD_p(\Gamma)$. To complete the proof of the theorem we must establish a bijection between $BHD_p(\Gamma)$ and $BHD_p(M)$, which we will now proceed to do.

Let $\bar{h} \in BHD_p(\Gamma)$. Define a function $h \in BD_p(M)$ by (\ref{graphtorm}). Denote by $\pi(h)$ the unique element of $BHD_p(M)$ given by the Royden's decomposition of $h$. Define $\Psi \colon BHD_p(\Gamma) \rightarrow BHD_p(M)$ by $\Psi(\bar{h}) = \pi(h)$. 
We will now show that $\Psi$ is one-to-one. Let $x \in \partial_p(\Gamma)$ and let $(x_n)$ be a sequence in $\Gamma$ that converges to $x$. Let $\Phi$ be the homeomorphism from $\partial_p(\Gamma)$ to $\partial_p(M)$ given in earlier in this section. Since $\Psi(\bar{h}) (\Phi(x)) = h (\Phi(x))$ for all $x \in \partial_p(\Gamma)$, it follows that 
\begin{equation} 
   \vert \Psi(\bar{h})(x_n) - h(x_n) \vert \rightarrow 0 \mbox{ as } n \rightarrow \infty. \label{hconvbd}
\end{equation}
Combining this with Lemma \ref{convinftyrm} we obtain $\vert \Psi(\bar{h})(x_n) - \bar{h}(x_n)\vert \rightarrow 0$ as $n \rightarrow \infty$. Thus $\Psi (\bar{h})(\Phi(x)) = \bar{h}(x)$. Let $\bar{h}_1, \bar{h}_2 \in BHD_p(\Gamma)$ and assume $\Psi(\bar{h}_1) = \Psi(\bar{h}_2)$. Then $\bar{h}_1(x) = \bar{h}_2(x)$ for all $x \in \partial_p(\Gamma)$. Hence $\Psi$ is one-to-one because a $p$-harmonic function on $\Gamma$ is determined by its values on $\partial_p(\Gamma)$. 

All that is left to do is show that $\Psi$ is onto. Let $u \in BHD_p(M)$ and define $u^{\ast} \in BD_p(\Gamma)$ by (\ref{rmtograph}). Denote by $\bar{h}$ the element in $BHD_p(\Gamma)$ given by the Royden decomposition or $u^{\ast}$. Let $y \in \partial_p(M)$ and let $x \in \partial_p(\Gamma)$ such that $\Phi(x) = y$. Pick a sequence $(x_n)$ in $\Gamma$ for which $(x_n) \rightarrow x$. Now $\Psi(\bar{h})(y) = \pi(h)(y) = h(y)$ since $y \in \partial_p(M)$. By (\ref{hconvbd}) and Lemma \ref{convinftyrm} we see that 
\begin{equation}
\vert \pi(h)(x_n) - \bar{h}(x_n) \vert \rightarrow 0 \mbox{ as } n \rightarrow \infty. \label{hnetconv}
\end{equation}
By Lemma \ref{infinityconv} 
\begin{equation} 
\vert u^{\ast} (x_n) - u(x_n) \vert \rightarrow 0 \mbox{ as } n \rightarrow \infty. \label{convu}
\end{equation}
It is also true that 
\begin{equation}
\vert \bar{h}(x_n) - u^{\ast}(x_n) \vert \rightarrow 0 \mbox{ as } n \rightarrow \infty, \label{huconv}
\end{equation}
due to $u^{\ast} = \bar{h}$ on $\partial_p(\Gamma)$. Combining (\ref{convu}) and (\ref{huconv}) with (\ref{hnetconv}) we obtain 
\begin{equation} 
\vert u(x_n) - \pi(h)(x_n) \vert \rightarrow 0 \mbox{ as } n \rightarrow \infty.
\end{equation} 
Thus $u(y) = \Psi(\bar{h})(y)$ for all $y \in \partial_p(M)$. Hence $\Psi(\bar{h}) = u$ by \cite[Lemma1]{Lee}. The proof of Theorem \ref{bijectionpharm} is finished.

\bibliographystyle{plain}
\bibliography{graphsmanpbd_rev}
\end{document}